# Werner DePauli-Schimanovich
# Institute for Information Science, Dept. DB&AI, TU-Vienna, Austria
Werner.DePauli@gmail.com

## The Notion "Pathology" in Set Theory (Paper), K50-Set5b

> *"Sets are the closed hereditary*
> *non-patho predicate extensions!"*[1]

### 1 Abstract
When we study the paradoxes of set theory we find out that there are mainly 2 types: the pathologies and the antinomies. These 2 notions are made precise and compared with the somehow inductively definable concept "abnormal". (See K50-Set7.) In the following 5 Patho Theses are discussed in order to formalize this notion of pathology. This allows us to define formally the property "Hereditary-non-Pathological" for well-formed formulas. With this property the system NACT∗ of Naïve Axiomatic Class Theory is constructed, which has a "unique maximal" universe (in a special sense). The references of all set-theory papers in this book are at the end in K50-Set10.

### 2 Introduction
Every set theoretician knows that most of the paradoxes of Naive Set Theory follow already directly from pure logic. But for several reasons nobody used this fact as decision criterion for the constitution of sets until now. (See K50-Set4.)

Of course there are also other antinomies like Cantor's or the complement of the ordinals. But these follow from proper axioms like separation or replacement. And if we forbid these axioms nothing bad can happen. Therefore we can define (for logical formulas) "Non-Pathology" this way and "Consistency" as "Hereditary-Non-Pathology" (i.e. non-patho for all sub-formulas[2]), which yields a new and interesting system of set theory.

The aim of this paper is to give a philosophical justification for the further use of (a restricted) Naive Mengenlehre (i.e. set theory) by the mathematicians. If an inconsistency can be derived from Naive Mengenlehre for a special "nearly-closed" formula A (i.e. a formula with the only free variable "x") as a substitution-instance for the infinite Comprehension Scheme, then it can already be derived from pure logic (augmented by the Comprehension Scheme with this substitution instance A and the axiom of extensionality).

---

[1] This is the most general principle for generating sets! The set-generating predicate has to be "nearly closed" or "almost closed", (which means here:) closed, except the variable x, which ranges over the objects satisfying the predicate.

[2] The suggestion to work out this concept in detail was made by Randall Holmes. He expects that this system would turn out to be equivalent to an existing axiom-system of positive set theory (if it is contradiction-free). In this case the condition "decidable", which I wanted first to use alternatively, can be dropped. If we want to exclude the undecidable cases, "Logo-cality" also can be defined in a similar way (as "Nearly-Closed Decidable Hereditary Consistency"). But it seems to be a misleading way to formalize Naïve Mengenlehre.

This fact gives rise to the following strategy: Of course we have to consider only nearly-closed formulas A (i.e. formulas without free parameters). But then it is enough to forbid the use of the "pathological" formulas as substitution instances into the Comprehension Scheme (i.e. those formulas which generate a contradiction already in the enriched pure logic with identity and function variables) and also all "sub-pathological" formulas (i.e. those which contain a patho formula as sub-formula in a special defined sense). This "hereditary-non-patho" property (if suitably worked out) may be already enough to generate a contradiction-free collection of set-constituting axioms.[3]

If this would not fit then we have to exclude also the undecidable formulas.[4] I.e.: The Comprehension Scheme is then restricted to (nearly-closed and decidable) hereditary-non-pathological formulas and this gives certainly a consistent system NAM* of "Naïve Axiomatic Mengenlehre"[5], but with less sets. If we formalize this thought within class theory we get the system NACT* of "Naïve Axiomatic Class Theory Star".

"Sets are the closed consistent (= CC) predicate extensions!" (Proper-classes are the inconsistent ones.) [Keep in mind: We use here the notion "consistent" as synonymous for "hereditary-non-patho".] This philosophical thesis will be explained in this paper below.

## 3 A Sequence of Systems
To reach this goal of a "consistent" (i.e. a hereditary non-patho) set theory we have to undertake a philosophical investigation about the question:

---

[3] This suggestion has been given to me by Randall Holmes in a private email in July 2002, because I was not sure that this would be the optimal solution. The exact definition is given in §10. In this case we have to forbid all other axioms, too (which are not substitution instances of the Comprehension Scheme of nearly-closed hereditary-non-pathological formulas), except extensionality and a (to small sets) restricted axiom of choice.

[4] If we do not demand decidability, it is still unknown if the CHnP-PE's (i.e. closed hereditary non-patho predicate-extensions) [or short: CC-PE's (i.e. closed consistent predicate-extensions)] yield a contradiction-free system. Therefore we mention decidability too, and call this system the CDC-PE's. But we do not study this system here in this paper.

[5] If we forget for a moment the "nearly-closed" condition for formulas (with only one free variable "x") because it can be lifted in principle, the hereditary-non-patho restriction is weaker then the stratification, which forms "New Foundations". (See Willard Van Orman Quine in the references.) But if we have to take (sub-)logo [= (hereditary-) non-patho & decidable] it becomes incomparable with strati, because there exist undecidable stratified wffs in NF too.

How can we find (step by step) a sequence of axiom systems for set theory (which converges to NACT*), where every "set-possible" class[6] or "as consistently conceivable" naïve set[7] can also be made a set in such an axiom system? This question is investigated in the following paper K50-Set6.

Our thesis (on which we hope to reach this goal) is to split up or divide the paradoxical sets (or predicates, or classes) in pathological, antinomical and abnormal ones. And the strongest philosophical thesis that will be formulated in this paper (of course always under the silent presupposition that the considered formula is nearly-closed) is that ONLY the (sub-)pathological formulas A are prevented from forming sets.[8] (If it would turn out that this restriction is not enough, we would have to forbid also the undecidable cases too.[9])

In case A is allowed to contain free parameters [in addition to the free variable "x"] one can construct some suitable wffs $A_j$ (j in J) such that a contradiction can be derived from their substitution instances into the comprehension scheme. In case A is Patho(A) you can already derive a contradiction from the Comprehension Scheme [with the substituted nearly-closed formula A having the only free variable "x"] by means of an extended pure logic. (In case A is undecidable, one of the formulas "A ==> B" or "non-A ==> B" such that Patho(B) will generate a contradiction.)

Following this philosophical line, "pathologicality" can be considered as a property of well-formed formulas (like stratification or closedness). Similarly to NAM* (i.e. the most important system of "Naïve Axiomatic Mengenlehre" NAM[10]), this will yield the system NACT* (i.e.

---

[6] A "set-possible" class or an "as consistently conceivable" naïve set is something like a parameter-free (and of course predicative) class that can be made to a set in some suitable system of set theory. (E.g. NF.) Later we can also allow free parameters if we can prove as meta-theorems that the classes generated in such way do not falsify the system.

[7] But the attempt of "Frege's Way Out" (which uses impredicative formulas to also render the Russell Class, Mirimanoff Class and Ordinal Class as sets) is excluded here. (See K50-Set8b.)

[8] There is still some additional work to be done to make this notion "(Sub-)Patho" and the corresponding part of logic (where e.g. functions should be allowed) precise. There is some hope to find for each A with Patho(A) a special decision procedure SDP(A). And since the existence of a general decision algorithm DA for this criterion Patho(A) turns out to be equivalent with the decision problem of the FOL, we are looking for DAs for large subsets of Patho(A). But as long we do not have such a DA, we work with Non-Monotonic Logic (what ZF and NF are in fact doing too).

[9] This causes some difficulties to define decidability for a nearly-closed wff A (i.e. with free variable x). My suggestion is: decidability of A(t) for all closed terms t. Because in most models there are more sets than closed terms, the addition of term-consistency (i.e. the corresponding rule to omega-consistency) would certainly be helpful.

[10] See: DePauli-Schimanovich, "NAM for Experiments" (K50-Set7b in this book) and "A Brief History of Future Set Theory" (K50-Set4 in this book).

the most important system of "Naïve Axiomatic Class Theory"[11]) if we use the Comprehension Scheme: "Every class {|x: A(x)|} with (nearly-closed & hereditary-)non-Patho(A) is a set". (During our first investigations non-Patho(A) should always mean Hereditary-non-Patho(A) and imply that A is a nearly-closed wff [i.e closed except only one free variable "x"].) And NACT* has no other additional axioms at all except extensionality and the axiom of choice for small sets!

All other systems studied in this paper or in K50-Set6 should only be viewed as a preliminary step in the sequence of axiom-systems to eventually reach NACT*, which we want to call colloquially "naked star". To construct this sequence, the notion of "Paradox" is split up into three different termini: "Pathology", "Antinomy" and "Abnormality", and a formal explication of these 3 termini is elaborated. We will start with the explication of the first terminus immediately and with the discussion of the other 2 termini directly after this.

## 4 Definition of Patho, Consi and Logo
Let us consider a formula-generator
FG == {A1, A2, . . . . Ai, A(i+1), . . . . . .}
which enumerates all nearly-closed well-formed formulas Ai (i.e. those [as decidable expected] wffs with "x" as the only free variable, later to be bound by the set-operator) of PIF-logic with the elementhood "in" as the only predicate constant and the set-operator "{x: . . . .}" as the only functional constant. (The PIF-logic is the FOL with identity and function symbols and reflexivity and substitutivity for the identity. Note that the set-operator is symbolized by the pure set-brackets, while the class-operator has a vertical bar in addition to the brackets.)

Furthermore we consider for all nearly-closed formulas Ai the relation
(CoSi) := forall y: [y in {x: Ai(x)} <==> Ai(y)].

(CoSi) is of course (for every Ai) a substitution instance of the general Comprehension Scheme
(CoS) := forall wff A: forall y: [ y in {x: A(x)} <==> A(y) ].

Every single axiom (CoSi) together with the axiom of extensionality (EE) and with PIF-logic we will call from now on CoSi-PIF. This is exactly the small logic fragment of set theory without any specific single axioms except two: (EE) and (CoSi).

If CoSi-PIF is inconsistent, we call the corresponding Ai pathological or patho for short. Once we have a "list" of all patho wffs (or a close super-set of them), we can easy define the terminus hereditary-non-patho in §10. If CoSi-PIF is non-patho and decidable we call Ai logo for short. Therefore:
for all i: Patho(Ai) or undecidable(Ai) or Logo(Ai).

---

[11] See: DePauli-Schimanovich, "NACT: A Solution for the Antinomies of Naïve Mengenlehre." (K50-Set4 in this book.)

E.g.: A1 := Verum and A2 := Falsum are logo. But we know of course that Russell-type formulas like "x non-in x" are patho. Our endeavor here is to find a formal characterization of the notion "Consistent" (i.e. hereditary-non-patho) or better "Logo". I.e. we are searching for the collections "Consi" (= HnP) of consistent formulas Ai (or "Logo" of all consistent and decidable formulas Ai) for which the following condition holds:
CoSi-PIF |--/-- (= non-derivable) Falsum.

"Consi" (or "Logo") shall also be used for the class of all sets formed with a consistent (or a logo) wff, and it is therefore equal to the universal set "us(NAM*)" of NAM* (and also NACT*). In analogy to the formulas Aj in Patho, the name "Patho" should also stand for the 2-nd order class of all proper-classes. (Note that the 2-nd order class-operator is symbolized by the set-brackets together with a slash, i.e. "{/" and "/}" !) I.e.:
Patho := {/ A: A is {|x: Aj(x)|} for a suitable j with Patho(Aj) /} .

Let us assume in the following the hypothesis that there exists a decision algorithm DA for all nearly-closed well-formed formulas A if they are patho (or logo).[12] (The real existence of such an algorithm is in fact irrelevant for our thesis.) Take an arbitrary nearly-closed formula A. Case Patho(A) will in general be easier to show (because a contradiction can easier be derived than the proof be given that no contradiction is possible). Also the case Decidable(A) is difficult and without special interest for us. Therefore let us concentrate on the case Non-Patho(A) and Herditary-non-Patho(A). DA will probably work in the following way: Since logic is consistent, a contradiction can only be derived by DA from (CoS-A). Every contradiction has the form "B & non-B" or "B <==> non-B". So either (CoS-A) produces the second prototype of contradiction directly or the first one follows from it.

Let e.g. A := "x in x"; if a contradiction of the 2-nd prototype is at all derivable, then only by the substitution of the complement[13] "ko(Ru)" of the Russell-class "Ru" into the variable y in (CoS-"x in x"). But nothing happens in this case, because this yields only the identity "B <==> B". For all other closed set-terms "t" substituted into y, either "t in t" is right or wrong, in which case t belongs to ko(Ru) or not. That property is either a fact of the system or can be chosen independently and cannot force a contradiction. (It is nothing other than the usual extension of predicates.) And the sets in ko(Ru) do not interfere with each other. Therefore nothing can happen, and also the first prototype of the contradiction cannot be deduced. In the other cases the argumentation will be similar, but of course more complex if we have complex predicates or formulas A. This is of course only a vague consideration.

---

[12] Since the class of all nearly-closed non-patho A is generated-up in an non-constructive way (similarly to the maximal model in the completeness proof of FOL) we do of course expect that an effective procedure does not exist. But this does not matter at all. It would be already enough to find a sufficiently large class of formulas for which a general DA can prove its hereditary non-pathologicality. But this does not influence our goal to investigate the consistency of NAM* and NACT*.

[13] Complement is always abbreviated with "ko" instead of "co", because the letter "c" is already reserved for class, comprehension, Church, choice, cofinality, consistent, closed, etc.

In these cases it would be very helpful if we also have to be sure that A is decidable. But since already the definition of the decidability of a nearly-closed wff Aj causes difficulties (because of the free variable "x"), we shift our considerations to the corresponding predicate extension (CoS-Aj) which is (now totally) closed after the substitution of Aj. And since with Aj(t) also non-Aj(t) is undecidable (for a special closed term "t" substituted for "x"), also (CoS-Aj) and (CoS-non-Aj) are too (and must be checked for all "t"). Therefore it is enough to investigate the decidability of (CoS-Aj) alone. And since we are already trying to derive a contradiction from (CoS-Aj), we shall do in parallel the same for non-(CoS-Aj). Of course, if Aj is undecideable, none of these 2 parallel machines will stop. But if non-(CoS-Aj) yields a contradiction, we know at least that Aj constitutes a set.

Independent of whether we have such a decision algorithm DA for the properties Consi(A) (or Logo(A)) of every (nearly-closed) well-formed formula A or at least for a large subclass of Consi or Logo, we can formulate our thesis: The general Comprehension Scheme (CoS) of Naïve Mengenlehre, restricted to only the (nearly-closed) A with Consi(A) (or Logo(A)), forms a consistent system of set theory, that we want to call "Naïve Axiomatische Mengenlehre" NAM*.[14] Let us define:
(Consi-CoS) := $U_i$ (CoS-i)   [= Union over i of the (CoS-i) with Consi(Ai) ].
Then our thesis is the conjecture:
Cons[PIF & (EE) & (Consi-CoS)]. (And the same should hold for Logo.)[15]

Or contentually in other words: "The contradictions of Naïve Mengenlehre can be derived already in pure logic." I.e. from a CoSi-PIF. Or as contra-position: If no contradiction can be derived from any CoSi-PIF for a suitable set F = {Ai: i in I} of nearly-closed formulas Ai, then the list of axioms "<(CoSi)> with i in I" forms a consistent system of Naive Mengenlehre. (In this case F =< Consi holds.)

There exist many systems of NAM. The most important one is the described system NAM* consisting of (Consi-CoS) and (EE) [and eventually an appropriate (AC) for small sets], but not containing any proper axioms like separation or replacement. Axioms (with free parameters) can only be shown as meta-theorems, e.g. the existence of the union or intersection or tuples of finitely many arbitrary closed sets. Or separation restricted to suitable small domains (like well-founded, Cantorian or "slim" sets) follows logically from (Consi-CoS), and probably restricted replacement follows too. Even the complement of "slim" and "mighty" sets can be shown to exist.[16] And that is more then enough for the mathematicians, because nobody wants to separate the Russell class "Ru" from the universal set "us".

---

[14] This means: Consi (and Logo) or better its subclasses are similar set-constituting formula-properties like "Stratified" of "New Foundations". The same holds for NACT*.

[15] NAM* (or NACT*) is the formalization of the method mathematicians work with in practice, already since the creation of Naïve Mengenlehre by Cantor. Antinomies appeared only for 2 reasons: 1-st because logicians like Bertrand Russell took into account not only the closed consistent predicate extensions (CC-PEs), and 2-nd because Cantor wanted to save separation as a general principle for arbitrary sets (and did not restrict it to only small sets). Since the working mathematicians did not produce any inconsistency we can consider this fact as 100 years soundness of NAM* (or NACT*) garanted by mathematical experience.

It is also probable that a suitable (AC), of course also restricted to small basic-sets, can consistently be added to NAM*, as it is done for NF or (unrestrictedly) for ZF.[17] We will call this enlarged NAM* then
NAMAC* := NAM* & (Restricted AC) [or for short also only NAM*, if the addition of the restricted (AC) does not cause any difficulties].

One can see that NAM* "contains" (in a certain sense) ZF (with restricted separation and replacement) (and the same is valid for NF), and for NAMAC* the corresponding fact holds.[18] (Consi-CoS) can also be formulated as "Every class {|x: A(x)|} with Consi(A) is a set." In similar cases we talk about "Naïve Axiomatic Class Theory" NACT and call this special system NACT* (= the naked star). In the "*"-case NAM* and NACT* coincide or yield the same system. But we will investigate other systems NACTi which are different from the corresponding systems NAMi. In the following we will concentrate on NACT.

## 5 Comprehension Schemes
Georg Cantor's Mengenlehre was mainly based on one principle: "Sets are Predicate-Extensions."[19] In formalized form this principle is known as the (inconsistent) Comprehension Scheme (CoS) of the Naive Set Theory for sets:
Forall wff A: forall y: (y in {x: A(x)} <==> A(y) ).
Note that the Russell-Paradox follows directly from this (CoS) without use of any other axioms, while other antinomies like Cantor's need the help of other axioms of set theory too.

---

[16] Definitions of "Slim" and "Mighty" are given later.

[17] See Kurt Gödel: "The Relative Consistency of the Axiom of Choice and the Generalized Continuum Hypothesis."

[18] "Contains" means the inclusion in a certain sense. Of course ZF generates "undecidable" sets too, but NAM* does not.

[19] We add only the 2 words "closed, consistent". (Or if it turns out that this will not be enough the 3 words "closed, decidable, consistent".): "Sets are [exactly the] closed consistent predicate extensions!" (Of course extensions of sets, because we have only this one type of entity!) To produce such sets we use the set-operator {x: A(x)} for the consi wffs A. But classes are arbitrary predicate extensions: "For all wffs A the class-operator {|x: A(x)|} constitutes a class. (Of course only of sets! Therefore we have 2 types of entities.)

After the discoveries of the paradoxes it was clear that restrictions have to be made to (CoS). The most prominent one is that of Willard Van Orman Quine: The formulas A in (CoS) have to be "Stratified", i.e. their variables should be indexed by natural numbers as in the Theory of Types.[20] This (StratCoS) is the basis for the system NF (= New Foundations). But for several reasons, the community of mathematical logicians has never accepted NF, despite the fact that (StratCoS) is a real Comprehension Scheme: it produces sets, and not only classes! And therefore NF is a restricted Naïve Mengenlehre.

Instead of using NF, mathematicians and logicians took another way: Every formula A in (CoS) should have been restricted to special sets. Let us call them for the moment the normal sets. Then we get the following scheme (NormCoS):
forall y: y in {x: A(x) & Normal(x)} <==> A(y) & Normal(y) .

But this is no longer a Comprehension Scheme, because it produces only the empty set 0, as long as we do not know which sets are normal. Therefore a shift was made in the systems and nomenclature: the old sets are, in the new system, a new entity called classes (symbolized by large Latin letters like X, Y, Z, .....), whereas normal sets are the new sets (for which we reserve now the small Latin letters x, y, z, ....), and the set-operator becomes the class-operator (restricted to sets). Class Theory was born by that act. But the old (CoS) is now only a Comprehension Scheme for classes, which are in fact virtual objects that can be eliminated in principle, if we use the scheme as something like a definition. (Or if we allow quantification over classes, then we have in fact a 2nd order logic.)

To differentiate this new scheme from the old (CoS) it has therefore been called the Church-Schema (CS):
Forall wff A: forall y: y in {|x: A(x)|} <==> A(y) .
Note that {|x: A(x)|} is now the class-operator forming classes, and the set-brackets are always together with a vertical bar. (CS) can also be formulated as:
Forall Y: Set(Y) ===> [Y in {|X: A(X) & Set(X)|} <==> A(Y)] ,
where Set(X) is a primitive predicate, which also could be replaced by "exist Z: X in Z", if quantification over classes were allowed.

The class-operator (C-O) together with the Church-Schema (CS), the extensionality axiom (EE), and maybe a suitable axiom of choice (AC), probably a restricted (to small sets) axiom (SmallAC), form the "frame" of class theory CT. But this frame CT does not build up any sets. For example the empty class 0 := {|x: Falsum|} can be equal to the universal class UC := {|x: Verum|}. To establish the existence of sets we need still the Zermelo-Fraenkel set-existence axioms. ZF together with the four axioms CT forms the Neumann-Bernays-Gödel class theory NBG.

---

[20] For the exact definition of stratification see [Forster 1995] and [Holmes 1998].

But CT can be used together with other set-existence axioms too. The systems of Naive Axiomatic Class Theory NACT, e.g., use completely different set-generators. The first system we want to investigate is NACT* which uses the (Consi-CoS) as Comprehension Scheme. But in this case we should keep in mind that, with the terminology-shift from sets to classes, the definition of subset (i.e. the subclass of a set in the case that it is a set) also changes too. It became the definition of a sub-class (of a class):
X =< Y := forall z: (z in X ==> z in Y).
This does not matter as long as we consider only "small"[21] set theories as ZFC or NACT#.

But if we allow "large" sets like a universal set too, the situation is different, because a subclass need not be a subset, e.g. Ru =< us (the Russell class is a subclass of the universal set), but ru =/< us (the class Ru is no subset, and therefore no set ru). The cases where a subclass is a set are regulated in the axiom of separation, which is usually restricted to small sets in universal set theories. If we want to express the subset relation, we need therefore small letters (on the left side):
x =< Y := forall z: (z in x ==> z in Y) := Set(X) & X =< Y .

A similar situation happens with the power-class2 (which can be a class of 2nd order in the general case), the (normal) power-class (which can be a proper-class too), and the power set (of a class or a set):
PC(Y) := {/ X: X =< Y /}, while
P(Y) := {|x: x =< Y|}, or
P(y) if Y is a set, or
p(Y) := {x: x =< Y} if P(Y) always constitutes a set, and
p(y) := {x: x =< y}, the same as the previous formula, but with the set y instead of the class Y. The small p (or the word "power") should indicate that the result of the power-operation is always a set.

Having this distinction in mind, it is evident that Cantor's antinomy cannot be derived in NACT* [despite the fact that P(us) = us], because every ko(Di(f)) or ko(Di-closed) of P(us) is no set.[22] Before we can define NACT* exactly, we have first to undertake a philosophical analysis of the terminus "Paradox".

## 6 The Term Paradox and the first 2 Patho-Thesises
In the following we want to split the terminus "Paradox" and consider the 3 termini "Pathological", "Antinomical" and "Abnormal". These termini are completely different. Pathological(A) means that the formula A (or rather the assumption that "the class generated by it is a set") leads to a contradiction independent from any axioms of a set theoretical system (except the extensionality axiom) only on the basis of (CoS) and pure logic. Pathological is therefore something like an "absolute" notion in the sense of Gödel.

---

[21] "Small" may be Founded or Hereditary-Founded, Mirimanoff, Cantorian, or Slim (i.e. not equal mighty to its complement). Concerning NACT# see "Naïve Axiomatic Class Theory: A Solution for the Antinomies of Naïve Mengenlehre" in this book K50-Set6.

[22] The exact definitions of these notions of diagonalization are given later in chapter 9.

We understand immediately that Russell's paradox is in fact a pathology. Also Mirimanoff's paradox, Burali-Forti's, etc, and also all formulas constructed in a similarly special way -- a point discussed later. On the other hand, Cantor's paradox is not pathological at all, because it can only be derived by use of the axiom of separation (or better replacement). Also the fact that the complement ko(On) of the Ordinal-Class On is again a proper-class can only be derived by the use of NBG-axioms. These paradoxes we will define later as antinomies.

We will use later for some axioms implicitly the Pathology Thesis. It splits into several points. The 1-st one is (1PT):
"Of every pair "A" and "non-A" (or the corresponding classes, generated by it) are not both pathological!"

We can formalize this 1-st thesis and call it (1JT). But first we have to give a formal Definition (DefPatho)[23] of Patho:
Patho(A) :<===> non-forall beta: [beta in {alpha: A(alpha)} <==> A(beta)] ,
and Patho({|x: A(alpha/x)|}) :<===> Patho(A) .
If we can derive the righthand side property for A in CoS-A-PIF, then we have shown Patho(A). Then (1JT) can be formulated in these terms:
non-Patho({|x: A(x)|} ) or non-Patho({|x: non-A(x)|} ).

The reader should keep in mind that "A" and "non-A" are here nearly-closed formulas without free parameters (except the "x"). As usual we mean with Patho also the 2-nd-order class of all Patho classes.[24] The (decidable) negation of Patho we call Logo. With these notion the 1-st Patho-Thesis becomes very short:
[Decidable(A) ==>] Logo(A) or Logo(non-A) .

We took here the greek letters because nobody should mix up the Comprehension Scheme (CoS) of Naïve Mengenlehre with the Church-Schema (CS) of the class-theoretical frame CT. (DefPatho) is no formula of NBG or NACT*, but of first-order logic with equality (and maybe additional function symbols), predicate constant "in", extensionality axiom, and class-operator (that may even be an "eliminable" notion). (DefPatho) means: if the comprehension scheme of naïve set theory delivers a contradiction only by use of an extended pure logic, then the set-generating predicate is pathological and also the class generated by it.

---

[23] We use here small Greek letters to remind the reader that these are naive sets of the (CoS) and not classes! Note also that the operator's parenthesis (in the 1st line of the definition) do not have an additional vertical bar (which means that this is the set-operator in contrast to the class-operator of the 2nd line of the definition).

[24] Consi (i.e. Herditary-non-Patho) is the 2-nd-order class of all consi classes, which will turn out to be a 1-st-order proper-class (or even a set) if Consi is contradiction-free.

For other axioms we will assume the 2-nd Patho-Thesis (2PT):
"The slim and mighty classes are non-patho." Formalized that yields (2JT):
Slim({|x: A(x)|}) or Mighty({|x: A(x)|}) ===> non-Patho({|x: A(x)|}).

We use here the following abbreviations:
Slim(X) := card(X) < card(ko(X)) and
Mighty(X) := card(X) > card(ko(X)) , where
Card(X) is the cardinality or cardinal number of X, or: card (X) if X is a set.
Ko(X) := {|y: y non-in X|} is the complement of X, or: ko(X) if Ko(X) is a set.
Slim := {|x: Slim(x)|} is the class of all slim sets and the corresponding scheme is used for "mighty". It is also possible to use instead of Slim the notion Cantorian of NF:
Cantorian(X) := card(X) < card(power(X)), where power(X) means the power set od X. (Of course in NF small latin letters are used for classes or stratified sets. In NF the power-property is stratified and the power-class therefore a set.)

As set-theoretical characterization I conjecture that Patho is the collection of all classes which cannot be a set produced as complement in any suitable system of set theory. [This follows from (1JT) and (2JT)]. (There are of course systems where a properclass of NBG can be a set [e.g. in K50-Set8b], but not produced as complement. The intended meaning of this conjecture is described in more detail in "Naïve Axiomatic Class Theory", K50-Set6, Chapter 7 in this book.

## 7 The notions Antinomic and Abnormal

The next type of paradox is "Antinomic(A) according to an axiom-system S", and it means that the assumption that the class generated by A is a set leads to a contradiction with the axioms of S. E.g. the predicate "Verum" [for the assumption of a universal set] is antinomic according to ZFC, or the formula "Ord(x)" according to NBG[25] [since we know that in NBG also "non-Ord(x)" produces a proper-class when we use axiom of replacement]. Thus we can see that an analog "Antinomy Thesis" to (1PT) or to (1JT) would be wrong.

But in analogy to (DefPatho) we can give a meta-definition for "Antinomic(X) with respect to system S" and call it (DefAntinom-S):
Antinom(sub S)({|x: A(a/x)|}): <===> non-Patho({a: A(a)}) & [proper-class({|x: A(a/x)|}) ==> Cons(S)] & [(assumption) Set({|x: A(a/x)|}) ==> non-Cons(S)].

Here again "a" should remind us of alpha and be a set of naïve set theory, and {a: ...} (without the bars before and after the set-brackets) is the naïve set-operator. It can even happen that the derivation of Falsum is accomplished by showing Patho({y: B(y)}) for some other class assumed to be a set. The negation of Antinomic shall be called Nomic.

---

[25] Note that the contradiction depends on the way the ordinals are defined! In NF the ordinals NO form a set.

Antinomic predicates are of no interest for us, because for every system S for which the property A is Antinomic there exists another system T for which it is not. E.g. for "Verum", which produces Cantor's antinomy of the universal set in NBG (if we assume that the universal class UC is a set), there exist many of systems where "Verum" generates a set.[26] If you restrict the axiom of replacement in ZFC (and of course the axiom of separation, which can be derived from it), you can also add complementation to ZFC, yielding a system ZFCK with a universal set.[27] Therefore the "Limitation of Size"[28] is not a feasible intuition in our case, because with another intuition such proper-classes as the universal class can be rendered into to sets. And our goal for NACT* is to render as many classes as possible into sets. (If they do not generate a contradiction already in pure logic, then they will also not generate a contradiction in NACT*.)

The third paradoxical terminus technicus is "Abnormal(A) with respect to the system S", and it is a property just to be determined implicitly by the system S. Therefore the primitive Predicate "Normal(x)" (which is of course the negation of Abnormal) has to be added to the language of the system S. E.g. "A = Verum" makes the universal set in the system NAM# abnormal in spite of the fact that Verum is not pathological at all. The predicate "Normal" is used in systems of NAM (= Naïve Axiomatic Mengenlehre) and can be eliminated from the language in principle. (We use it only for didactical reasons here, because it makes the working of the systems better understandable.) An analog meta-definition for "Abnormal(A) respectivly to system S" can be given if we replace (CoS) by the "restricted (CoS) of S", i.e. of the restricted Comprehension Schema of some system of NAM.

It would lead us too much into technical details if we want to investigate here in this paper the different systems of Naïve Axiomatic Mengenlehre,[29] too. But we want to give the reader 2 examples for such a restricted (CoS) to see the comparison with NAM (which is only the same as NACT in the star-case). The first system we mention is NAM% with:
(raBaDi-CoS): = forall y [y in {x: A(x)}  <==>  A(y) or non-Normal ({x: A(x)})].

---

[26] Thomas Forster's excellent book gives a detailed overview over these systems. [Forster 1995].

[27] See http://info.tuwien.ac.at/goldstern/papers/notes/zfpk.pdf
The paper of Martin Goldstern contains a proof for a slightly weaker conjecture than the one above.

[28] See the book by Hallett [1984]: "Cantorian Set Theory and Limitation of Size".

[29] See DePauli-Schimanovich [2006]: NAM for Experiments. (See K50-Set7b.)

NAM% together with (NSxorK):= [Normal(x) <=/=> Normal(ko(x))][30] and (NC2) := (SPK-Cond) := [Slim(x) ==> Normal(x)] yields NAM# (and with the 4 fundamental axioms NAM#4).
The second system is NAM& with
(rinoBaCo-CoS) := forall y [y in {x: A(x)} <==> A(y) & Normal ({x:A(x)})]
The third system is NAM§ with:
(NorBI-CoS) := forall y (Normal({x: A(x)}) ==> [y in {x: A(x)} <==> A(y)]) .

NAM§ together with (NSoK):= [Normal(x) or Normal(ko(x))] and (NC6) := (nEMK-Cond) := [x ~/~ ko(x) ==> Normal(x)][31] yields NAM+ (and with (StratNor) the system NAM+Strat). Instead of (StratNor) we can use also (CondNor): = [Condition(x) => Normal(x)] .

Condition(x) can be any suitable condition, e.g. "x =/= @", which matches all abnormal sets with the commercial at "@". (In this case we need additional set-existence axioms.) If Condition := Logo({x: A(x)} ), NAM§ becomes NAM*.
We can see that in such cases the terminus Abnormal is defined implicitly by the system of NAM involved.[32]

## 8 Hereditary-non-Pathological Predicates
The goal of this paper is to establish a set theoretical system NACT* as an extension of the class theoretical frame CT with the Comprehension Scheme
(ConsiCoS) := [Hereditary-non-Patho(A) ==> Set( {|x: A(x)|} ) ].
Only if Consi turns out to be not strong enough we use:
(LogoCoS) := [Decidable(A) & Hereditary-non-Patho(A) ==> Set( {|x: A(x)|} ) ].

The only additional axioms of NACT* are the axiom of extensionality (EE) and an axiom of choice (SlimAC) restricted to "slim" sets.

A first approximation to the notion "non-Patho" is the stratification of formulas. Stratification is an easily decidable property of a given formula, non-Patho is not. Probably it will be possible to find an algorithm which can reduce all (nearly-closed) formulas A with Patho(A) to Falsum and to show for all A with non-Patho(A) that no contradiction is derivable. Than also all A with Consistent(A) are characterizable. (But it is not sure that such an algorithm exists also for the A with Logo(A), at least for an appropriate subset.)[33]

---

[30] Class-theoretically this reads (SetDisSetInko):=[Set ({|x: A(x)|}) or Set ({|x: non-A(x)|} ) ], but in this case not both.

[31] In the language of class-theory this is (2JT).

[32] I have investigated about 100 different systems of NAM.

[33] I want to work with 2 PhD-students to try to find a suitable large class of predicates which solves this decision problem for Patho.

As long as we do not have any solution for this decision problem for Patho, I suggest working with Non-Monotonic Logic. In fact ZF and NF also do this, assuming their axioms are valid and facing for a possible contradiction, and in this case we repair the inconsistency, as Quine did with ML. Since Gödel's incompleteness results[34] we have no other possibility either. Furthermore, this is exactly the way mathematicians work in practice. They do not construct sets with ZFC but still use Cantor's Naïve Mengenlehre. If their intuition tells them that a special predicate is not pathological, they form a set with it, and a contradiction never occurs. This working method of the mathematicians is somehow an experimental proof or a confirmation in practice of the correctness of NACT*.

Since (for nearly-closed A) "Stratified (A) ==> Hereditary-non-Patho (A)" is evident, "(StratCoS) =< (ConsiCoS)" follows, because (EE) is the same in both systems. (That this holds also for parameterized A's has to be shown by meta-theorems.) Therefore: "NACT* ==> NF$^{()}$" is probably valid too, even with the axiom of counting. (Of course for (LogoCoS) this is only valid for decidable statements.)

You can also add a restricted axiom of choice to both systems, because NF is compatible with (CantorianAC) where the basic set must be Cantorian. But it is unknown whether this implication holds when you are adding the axiom of counting to NF. Therefore let us start with stratification as a first approximation to hereditary-non-pathology and let us look for possibilities to refine it.

The first question is: what happens with the unstratified non-pathological formulas and their classes? Of course according to axioms like (SetDisSetInko), at least one of a positive and negative predicate must generate a set. But (in contrast to NACT# or NAM#) we want this time that both generate a set, if possible. There are unstratified predicates where both formulas, "A" and "non-A", can generate a set consistently. Implicitly Quine's Thesis is that if a formula is stratified, then it is not pathological (or antinomic). But we strengthened this thesis and formulated the "Hereditary-non-Patho (is a) Set" axiom (HereNonPathoSet):
(ConsiCoS) : = Closed Hereditary-non-Patho({|x: A(x)|}) ==> Set({|x: A(x)|}).

## 9 The third Patho-Thesis
The axiom schemes (ConsiCoS) (used as our only axiom!) together with the CT-frame forms the optimal system NACT*. (The asterisk "*" refers to the optimal solution, as in game theory.) It renders every class, which need not be a proper-class to a set. Therefore NACT* generates a unique maximal[35] universe, which is furthermore the maximal universe in an absolute sense. (If A is very complicated, it may be that we cannot decide ad hoc whether Patho(A) holds or not. In this case we work with Non-Monotonic Logic and assume non-Patho(A), as long a contradiction do not arise). And since we have no other axioms than the CT-frame and (ConsiCoS), no antinomy can arise (because we defined the antinomies as those contradictions which are not already a pathology).

---

[34] See all my Gödel-books and the film in the literature.
[35] Maximal should mean that if you add a class (which satisfies the axioms of CT and) which is not yet in the universe as a set, the system becomes inconsistent.

Axiom scheme (ConsiCoS) presupposes of course the 3-rd Patho-Thesis (3PT) that "The Hereditary-non-Patho classes can coexist together as sets, not influencing each other and therefore not yielding an inconsistency for a system S. I.e. an inconsistency of a class-theoretical system S can only arise from the axioms of the system and from pathological classes. If we have no axioms except (ConsiCos) as in NACT*, the system must be consistent."[36] Or: "If X is a class of non-patho sets y, living consistently together in peace, we can add any hereditary-non-patho set z, such that "X u {z}" will not become inconsistent. This means: "We can add all closed hereditary-non-patho classes step by step as sets to our present universe of NACT* (starting with 0) by transfinite induction, until it becomes complete and contains all closed hereditary-non-patho sets."

Formalized, this 3-rd Patho-Thesis is (3JT), where UCnascendi(NACT*) is UC considered to be in status nascendi:
[forall y: y in X ==> Closed & Here-non-Patho(y)] & X =< UCnascendi(NACT*) & Cons(NACT*(UCnascendi)) ======> forall z: [Closed & Here-non-Patho(z) ==> Cons(NACT*(UCnascendi u {z}))]. Or:

forall y: [y in UCnascendi(NACT*) ===> Closed & Here-non-Patho(y)] ===> Cons(NACT*).

From this it follows that the implication in the brackets can be strengthened to an equivalence, when the status nascendi is finished and the UC is complete:
UC(NACT*) = Consi ====> Cons(NACT*).

If we investigate the complement ko(Di(f)) := {|x: x non-in f(x)|} of the "diagonal" class Di(f) in Cantor's proof of "card(x) < card(power(x))", we can see that the formula generating ko(Di(f)) is predicative. But the proof makes the function f impredicative. This is not forbidden in ZFC, NBG or the NACT-systems (as long as no inconsistencies can be derived from this impredicativity), but it causes some difficulties for NACT*. For a fixed x0 of the domain of f, f(x0) shall be (in Cantor's proof) equal to ko(Di(f)). Since in the (DefPatho) the class-generating predicates should be nearly-closed formulas, we have to consider:
ko(Di-closed) := {|x: forall f: [x in domain(f) & f(x) =< domain(f) ==> x non-in f(x)]|}

If we put this class-term into the (CoS) for sets of Naïve Mengenlehre in (DefPatho), the right-side expression of (DefPatho) yields a contradiction, because we can move the general quantifier "forall f" outside together with the precondition (of the implication in the brackets).

---

[36] Of course this does not contradict Gödel's incompleteness theorem because the 3-rd Patho-Thesis (3PT) is no formal proof of the inconsistency of NACT*, but only a contentual argument. Concerning Gödel, consult the literature in the references.

Then we have to replace all occurences of beta by the x0, for which $f(x0) = \{|x: x \text{ non-in } f(x)|\}$ holds, which causes the contradiction in the rest of the formula (what remained) in (CoS). From this fact follows that no function from domain(f) onto power(domain(f)) exists. This means that Cantor's theorem is provable in NACT*.

Our problem is the question: Is ko(Di-closed) pathological? Concerning our definition of (DefPatho) it is not, because we moved the function quantifier and precondition outside. But therefore the derivations are similar to the ones of Cantor's proof, where A in ko(Di(f)) is not a nearly-closed formula. This is also a reason why I ask under what circumstances can we drop the condition in (DefPatho) that A has to be a nearly-closed formula?[37] But before we can decide this question, more investigation has to be undertaken into the study of similar term-constructions like ko(Di-closed).

In future we will drop the word "closed" in connection with HnP sets and assume it implicitly.

## 10 Circumscription of Patho within pure logic and the 4$^{th}$ and 5$^{th}$ Patho-Thesis

Let us try now to circumscribe the intensional meaning of Patho(A). First we know that "x non-in x" is patho. We want to call this proposition the "non-circular formula of step 1" or NC1(A). Then

$NCn(A) := [\text{non-exist } x0, ...., x(n>0): x = x0 \text{ in } x1 \text{ in } ..... \text{ in } xn = x]$ or non $(x \text{ in}^n x)$, where "in$^n$" means a n-steps-long epsilon chain. $NC_n(x)$ means the same formula with variable "x".

The NCn are the basic patho wffs. (I am convinced that nearly all other patho wffs are derivable from them.) Let us call e.g. a wff A self-applicable if $A(\{x: A(x)\})$ holds. [E.g. A := Set(X) if a universal set exists.] If such a wff A implies NC1, i.e. "$\{x: A(x)\}$ non-in $\{x: A(x)\}$", A is of course patho too. This is also the case for e.g. Miri(X)[38], Ordinalnumber(X), Cardinalnumber(X), etc. A similar situation happens if A implies NCn and A is "sub-self-applicable", i.e. for some chain $x = x1 \text{ in } x2 \text{ in } ... xn = \{z: A(z)\}$ the statement A(x) holds. We call all these A's (for the first) "derivatives of NCn".

Let us now consider NCn(X) and its derivatives like Miri(X), On(X), etc, i.e. the classes generated by the corresponding predicates are of course patho and proper-classes. Since "NCn(X) ==> Patho(NCn(X-{x}))" and the same with Miri and the other drivatives, all "NCn(x) & non-A(x)" are patho, where A(x) creates only a small number of exceptions. A similar characterization holds for "NCn(x) or A(x)", where A(x) is as before. The same holds for the Miri(A) and the other derivatives. This means you can add to or take out of a properclass a suitable number of sets of a class, especially if this class forms a set.

---

[37] As already mentioned e.g. union and intersection allow free variables, also the complement of small sets, etc. But we are only allowed to establish the existence of Hereditary-non-patho sets with free parameters by meta-theorems.

[38] Miri(A): = A is a formula saying that there exists no infinite descending sequence of elements. Note that we can express these features by pure logic (only by extending the language using set symbols).

Since our 1-st Patho-Thesis is "Patho(A) ==> non-Patho(non-A)" (where with A we make also a proposition about the "patolog[icalit]y" of the Negation of A), we can also ask what happens with the dualization D(A) of a formula A, or the (Contra-)Valuation or "Value-Change" V(A) of A (where every "in" in A is replaced by "non-in"). In what cases is D(A) patho together with A, and V(A) non-patho with non-A? This means that even the case with only small exceptions is not yet solved at all.

But we want to formulate the 4-th Patho-Thesis (4PT): "If the exceptions A of NCn or Miri are of less power than their own (i.e NCn or Miri), then Patho(NCn&non-A) and Patho(Miri&non-A) and also Patho(NCn or A) and Patho(Miri or A). And the corresponding fact holds for the other derivatives of NCn too." Formalized, this yields (4JT) where Deriv(NCn) shall denote a representative of all derivatives of NCn, including NCn itself:
card(A) < card(Deriv(NCn)) ===> Patho(non-A & Deriv(NCn), A or Deriv(NCn)).[39]

This 4$^{th}$ Patho-Thesis can be strengthen under the assumption that all proper-classes have the same cardinality and are furthermore equal-mighty to its complement, while sets have a smaller or larger cardinality than its complement. (Cf the (2JT).) This yields (5PT) and its formalization (5JT):
Card (A) < card (non-A) ==> Patho (non-A & Deriv(NCn), A or Deriv (NCn)).
It may be the case that if y on add a mighty set that the pathology disappears.

Let us call the class-theoretical frame CT plus (5JT) the system NACT$ where Patho(A) is always compatible with (5JT). But the real difficulties start when the number of exceptions (i.e. the class {|x: A(x)|} ) ~ Miri. ("~" = "of equal power as"). Miri can be split into several parts which are all patho. But what parts are so and what are not so (and can therefore be used as exceptions), is the great question here. That means, presently we cannot find a complete suitable (syntactical) characterization of Patho(A) (which is defined as derivability of Falsum from (CoS-A-PIF)) by pure syntactical structural properties of A. (Until now we need still some semantic properties of the classes generated). But some future geniuses will find a 6-th Patho-Thesis (6PT) and its pure syntactical formalization (6JT) as a solution.[40]

---

[39] This is a somehow circular statement because the cardinality of "Deriv(NCn)" in the antecedence of the implication depends on the property "non-Patho" and therefore "Patho" too. Here, the thesis that "all proper-classes are of same cardinality" is used implicitly.

[40] In spite of these difficulties we want to try to use "NCn", "Miri" and the other derivatives of NCn for a syntactical definition of a large class of "Hereditary-non-Patho" formulas. My suggestion is to consider (as a first step) all exceptions in NCn and Miri also as Patho, even in the case that the "statement with the exception" would not be Patho any longer. (Therefore we put the defined termini in quotation marks.) Thus let us define:

## 11 Syntactical characterization of Here-non-Patho

A first suggestion to define a suitable subset SyntHnP of Here-non-Patho formulas A only syntactically should be suggested here:

First we have to find the primitive pathologies Prim-Patho; these are of course NCn and Miri. But we know that Ordinal-Number On(x), Cardinal-Number Cn(x), Wellfoundedness Fund (x) and other properties of sets can also be Patho. So our first partial goal will be to find all Prim-Pathos. But our hope is that by the following definition they all can be reduced as derivates of NCn or Miri.

The idea is to forbid all formulas A, where some other arbitrary formulas $B_1 \ldots B_k$ are added at some suitable places in a NCn (for some n and in correspondence with the usual formula construction algorithm). The same should be valid for Miri. This way all expections "& non-Bi" and "or Bj" are produced.

So lets first fix a number $n \geq 0$ and the n variables x1 … xn. Then consider the formulas "x1 in x2", "x2 in x3", … "x(n-1) in xn".
NCn is then non exist x1, … xn: x in x1 & x1 in x2 & … xn in x.

(Def-SyntP):
A derivative of NCn should be a formula A where one or more formulas c of the form " & non-Bi" or "or Bj" are added to NCn:
A = non exist x1, … xn: x in x1C & x1 in x2C & … xn in xC.

All these derivatives A form the collection SyntP and should be included in the collection Deriv(NCn). The same has to be formulated with Miri. We hope that also On, Cn, Fund, etc. will be in Deriv(NCn u Miri). If we can show that this is the case, we are at the end. If not, we have to find all exceptions and formulate an analogue definition including these exceptions.

(Def-SyntHnP):
This is the collection of all wffs minus the syntactical pathos:
SyntHnP : = wff \ SyntP.

---

(1) Every nearly-closed wff A is "primitive-Patho" if it is constructed from NCn or Miri or any other derivative of NCn by addition of some wffs B1 to Bk (in the recursive hierarchical construction of the wff A) to NCn or its derivatives.
(2) The case of quantifying such a wff (composed of A1 to Aj) has to be treated concerning the rules of construction of primitive-patho wffs (already in the 1-st step).
(3) Such primitive-Patho formulas A can also be in a term (like the set operator) or logical particle (like the "ι" [this x such that A]).
(4) Now we need not any longer to give a proof that (1) is equal to Patho(A) [i.e. (CoS-A-PIF) yields a contradiction].
(5) A is "primitive-non-patho" if it is a (nearly-closed) primitive wff of FOL (augmented with "in") which is not "primitive-Patho". With this concept we can define the notion hereditary-non-patho [or short: consi].
(6) A is primitive-consi if it is primitive-non-patho.
(7) If A and B are primitive-consi wffs then non-A, A&B, etc are it too. (The only free variable is "x".)

The definition (Def-SyntP) is in fact an enlargement of (4JT) or (5JT). Its not longer demanded that card(A) of an exception-class {|x:A(x)|} is of smaller cardinality than NCn (or some of its derivatives) or than its complement. We exclude all exceptions, also if its cardinality is the same as NCn or Miri or a derivative of them. That makes (Def-SyntP) much easier and allows us the move away from semantical considerations to syntactical ones. Therefore also (Def-SyntHnP) becomes very easy.

It does not matter if by the stronger (Def-SyntP) also cases are excluded which would be allowed by (4JT). I.e. if card(A) ~ card(NCn) happens it is still possible that non-Here-Patho (NCn & non-A). We can always assume that some equivalent formula B exists which has nothing to do with a Deriv(NCn). B would us then allow to construct {|x:B(x) |} with the same extension as {|x:NCn(x) & non-A(x) |}. That means that B is in SyntHnP.

With (Def-SyntHnP) we define probably only a subclass (SyntCoS) of (ConsiCoS). But it seems that this subclass will be enough for our considerations to establish a maximal universe of sets in Naïve Mengenlehre. The syntactical version of NACT$ should be therefore CT&(SyntCos), what we in future mean when we are speaking about NACT$.

After this syntactical analysis of the notion "Patho" we can go back to NACT$ (which contains no other special axioms then (5JT)) and add an axiom (ExceptionSet) which is a better explication of (5JT) and which somehow can be considered as a propabilistic simulation of the schema (ConsiCoS). We do that to randomize the decision procedure of Patho(A) instead of using the strict logically and not randomly formulated schema (ConsiCoS), which we do not use in NACT$. (ExceptionSet) shall be a criterion for the probability whether a formula A is patho or sub-patho or not. This will help us to sharpen our insight as to the likelihood of being pathological. Maybe one day we can also add an axiom (Deriv(NCn)PartSet) to NACT$, corresponding to (6JT) as a solution of the 6-th Patho-Thesis concerning what the exceptions of Deriv(NCn) are, i.e. the pathological sets with power equal to a special Deriv(NCn) like Miri.

If we add also "axiom" schema (ConsiCoS) to NACT$, we get this way the system NACT$*, which is a combination of both. Of course all axioms of any specific system of set theory should become redundant in NACT$* (which means that NACT$* is the same as NACT* if the added axioms (ExceptionSet) and tentatively (NCnPartSet) and (MiriPartSet) do not produce some inconsistencies). But since the decision algorithm DA of Patho(A) might require as much as 100 years to compute, an additional aid for this decision (and these are the axioms mentioned above) is always to the good.

## References for the Set-Theoretical Articles, K50-Set10

Bernays, Paul & Fraenkel, Abraham
[1968] Axiomatic Set Theory. Springer Verlag, Heidelberg, New York.
Brunner, Norbert & Felgner, Ulrich
[2002] Gödels Universum der konstruktiblen Mengen. In: [Buldt & al, 2002].
Buldt, Bernd & Köhler, Eckehard & Stöltzner, Michael & Weibel, Peter & Klein, Carsten & DePauli-Schimanovich-Göttig, Werner

[2002] Kurt Gödel: Wahrheit und Beweisbarkeit, Band 2: Kompendium zu Gödels Werk. oebv&htp, Wien.
Casti, John & DePauli, Werner
[2000] Gödel: A Life of Logic. Perseus Publishing, Cambridge (MA).
Davidson, Donald & Hintikka, Jaakko
[1969] Words and Objections. D. Reidel Publ. Comp., Dordrecht.
DePauli-Schimanovich, Werner (See also: Schimanovich [1971a], [1971b]) and [1981]).
[ca 1990] On Frege's True Way Out. Article K50-Set8b in this book.
[1998] Hegel und die Mengenlehre. Preprint at: http://www.univie.ac.at/bvi/europolis .
[2002a] Naïve Axiomatic Mengenlehre for Experiments. In: Proceedings of the HPLMC (= History and Philosophy of Logic, Mathematics and Computing), Sept. 2002 in Hagenberg/Linz (in honour of Bruno Buchberger).
[2002b] The Notion "Pathology" in Set Theora. In: Abstracts of the Conference HiPhiLoMaC (= History and Philosophy of Logic, Mathematics and Computing), Nov. 2002 in San Sebastian.
[2005a] Kurt Gödel und die Mathematische Logik (EUROPOLIS5). Trauner Verlag, Linz/Austria.
[2005b] Arrow's Paradox ist partial-consistent. In: DePauli-Schimanovich [2005a].
[2006a] A Brief History of Future Set Theory. Article K50-Set4 in this book.
[2006b] The Notion of "Pathology" in Set Theory. Article K50-Set5b in this book.
[2006c] Naïve Axiomatic Class Theory NACT: a Solution for the Antinomies of Naïve "Mengenlehre". Article K50-Set6 in this book.
[2006d] Naïve Axiomatic Mengenlehre for Experiments. Article K50-Set7b in this book.
DePauli-Schimanovich, Werner & Weibel, Peter
[1997] Kurt Gödel: Ein Mathematischer Mythos. Hölder-Pichler-Tempsky Verlag, Wien.
Feferman, Solomon & Dawson, John & Kleene, Stephen & Moore, Gregory & Solovay, Robert & Heijenoort, Jean van,
[1990] Kurt Gödel: Collected Works, Volume II. Oxford University Press, New York & Oxford.
Felgner, Ulrich
[1985] Mengenlehre: Wege Mathematischer Grundlagenforschung. Wissensch. Buchgesellschaft, Darmstadt.
[2002] Zur Geschichte des Mengenbegriffs. In: [Buldt & al, 2002] und Artikel K50-Set3 in this book.
Forster, Thomas
[1995] Set Theory with an Universal Set. Exploring an Untyped Universe. (2$^{nd}$ Edition.) Oxford Science Publ., Clarendon Press, Oxford.
Gödel, Kurt
[1938] The relative consistency of the axiom of choice and of the generalized continuum hypothesis. In: [Feferman & al, 1990].
Goldstern, Martin & Judah, Haim
[1995] The Incompleteness Phenomenon. A. K. Peters Ltd., Wellesley (MA).
Goldstern, Martin
[1998] Set Theory with Complements. http://info.tuwien.ac.at/goldstern/papers/notes/zfpk.pdf

Hallet, Michael
[1984] Cantorian Set Theory and Limitation of Size. Oxford Univ. Press, N.Y. & Oxford.
Halmos, Paul
[1960] Naive Set Theory.   Van Nostrand Company Inc., Princeton (NJ).
Holmes, Randall
[1998] Elementary Set Theory with a Universal Set.
Vol. 10 of the Cahiers du Centre de Logique. Academia-Bruylant, Louvain-la-Neuve (Belgium).
[2002] The inconsistency of double-extension set theory.
http://math.boisestate.edu/~holmes/holmes/doubleextension.ps
Jech, Thomas
[1974] Procedings of the Symposium in Pure Mathematics (1970), Vol. XIII, Part 2, AMS, Provicene R.I. .
Jensen, Ronald Björn
[1969] On the consistency of a slight (?) modification of Quine's New Foundation.
In: [Davidson & Hintikka, 1969].
Köhler, Eckehart & Weibel, Peter & Stöltzner, Michael & Buldt, Bernd & Klein, Carsten & DePauli-Schimanovich-Göttig, Werner
[2002] Kurt Gödel: Wahrheit und Beweisbarkeit, Band 1: Dokumente und historische Analysen. Hölder-Pichler-Tempsky Verlag, Wien.
Kolleritsch, Alfred & Waldorf, Günter
[1971] Manuskripte 33/'71 (Zeitschrift für Literatur und Kunst). Forum Stadtpark, A-8010 Graz, Austria.
Mathias, A.R.D.
[1992] The Ignorance of Bourbaki. In: The Mathematical Intelligencer 14 (No.3).
Quine, Willard Van Orman
[1969] Set Theory and its Logic. Belknap Press of Harvard University Press, Cambridge (MA).
Rubin, Jean & Rubin, Herman
[1978] Equivalents of the Axiom of Choice.   Springer, Heidelberg & New York
Schimanovich, Werner,
[1971a] Extension der Mengenlehre. Dissertation an der Universität Wien.
[1971b] Der Mengenbildungs-Prozess. In: [Kolleritsch & Waldorf, 1971].
[1981] The Formal Explication of the Concept of Antinomy. In: EUROPOLIS5, K44-Lo2a (and Lo2b). Or: Wittgenstein and his Impact on Contemporary Thought, Proceedings of the 2$^{nd}$ International Wittgenstein Symposium (29$^{th}$ of Aug. to 4$^{th}$ of Sept. 1977), Hölder-Pichler-Tempsky, Wien 1978.
Schimanovich-Galidescu, Maria-Elena
[2002] Princeton – Wien, 1946 – 1966. Gödels Briefe an seine Mutter. In: [Köhler & al, 2002].
Scott, Dana,
[1974] Axiomatizing Set Theory. In: Jech [1974] .
Suppes, Patrick,
[1960] Axiomatic Set Theory. D. Van Nostrand Company Inc., Princeton (NJ).

Weibel, Peter & Schimanovich, Werner
[1986] Kurt Gödel: Ein mathematischer Mythos. Film, 80 minutes, copyright ORF (= Austrian Television Network), ORF-shop, Würzburggasse 30, A-1130 Wien.
Wittgenstein, Ludwig,
[1956] Bemerkungen über die Grundlagen der Mathematik / Remarks on the Foundation of Mathematics. The M.I.T. Press, Cambridge MA and London.